\documentclass{amsart}
\usepackage{amssymb}
\usepackage{latexsym}
\date{\today}

\newcommand{\Z}{{\mathbb Z}}
\newcommand{\R}{{\mathbb R}}

\newcommand{\N}{{\mathbb N}}

\newcommand{\beq}{\begin{eqnarray}}
\newcommand{\eeq}{\end{eqnarray}}

\newtheorem{lemma}{Lemma}[section]
\newtheorem{theorem}[lemma]{Theorem}
\newtheorem{prop}[lemma]{Proposition}
\newtheorem{coro}[lemma]{Corollary}

\def\ind{\mathrm{ind}}
\def\A{\mathcal{A}_k}

\begin{document}
\title[AR Subshifts: Linear Recurrence, Powers, and Palindromes]{Arnoux-Rauzy Subshifts: Linear Recurrence, Powers, and Palindromes}
\author[D.\ Damanik, L.\ Q.\ Zamboni]{David Damanik$\, ^1$ and Luca Q.\ Zamboni$\, ^2$}
\thanks{D.\ D.\ was supported in part by NSF Grant No.~DMS--0227289}
\maketitle
\vspace{0.3cm}
\noindent
$^1$ Department of Mathematics 253--37, California Institute of Technology, Pasadena, CA 91125, USA\\[2mm]
$^2$ Department of Mathematics, University of North Texas, Denton, TX 76203, USA\\[3mm]
E-mail: \mbox{damanik@its.caltech.edu, luca@unt.edu}\\[3mm]

\begin{abstract}
We consider Arnoux-Rauzy subshifts $X$ and study various
combinatorial questions: When is $X$ linearly recurrent? What is
the maximal power occurring in $X$? What is the number of
palindromes of a given length occurring in $X$? We present
applications of our combinatorial results to the spectral theory
of discrete one-dimensional Schr\"odinger operators with
potentials given by Arnoux-Rauzy sequences.
\end{abstract}

%
%
%
%

\section{Introduction}

Mainly motivated by the discovery of quasicrystals by Shechtman et al.\ in 1984 \cite{sbgc}, there has been a lot of research done on the spectral properties of Schr\"odinger operators with potentials displaying long-range order. The first rigorous mathematical results were obtained in the late eighties. By now, many key issues are well understood, at least in one dimension. The two survey articles \cite{d2,s2} recount the history of this effort up to 1994 and 1999, respectively.

The primary example is given by a discrete one-dimensional Schr\"odinger operator whose potential is given by the Fibonacci sequence. More generally, one considers Sturmian potentials or potentials generated by (primitive) substitutions. It turned out that all these potentials lead to the same qualitative behavior: The corresponding Schr\"odinger operator has purely singular continuous zero-measure Cantor spectrum. This has been established for all Sturmian potentials and most substitution potentials. On the other hand, no counterexample is known. This led to the conjecture that these properties are shared by a large class of potentials displaying long-range order in a certain sense. One possible way to measure long-range order is given by the combinatorial complexity function $f : \N \rightarrow \N$ associated with a potential taking finitely many values, where $f(n)$ is given by the number subwords of the potential of a given length $n$. Since periodic potentials are well understood, one is interested in the case of an aperiodic potential and, in this case, it is well known that the complexity function grows at least linearly. One possible point of view could be the stipulation that long-range order manifests itself in a slowly (e.g., linearly) growing complexity function, possibly along with further conditions.

This combinatorial approach is further motivated by the fact that the properties above, singular continuous zero-measure Cantor spectrum, can be shown by purely combinatorial methods. The key combinatorial properties that allow one to deduce these spectral properties are linear recurrence and the occurrence of local symmetries such as powers and palindromes. Here, a sequence is linearly recurrent if its subwords occur infinitely often, with gap lengths bounded linearly in the length of the subword. Powers are repetitions of subwords and palindromes are subwords that are the same when read backwards.

Thus, the interplay between the spectral theory of Schr\"odinger operators and combinatorics of infinite words has enjoyed quite some popularity recently, due to its success in answering long-standing questions (e.g., the completion of the analysis of the Sturmian case \cite{dkl} or the proof of zero-measure spectrum for all primitive substitution potentials \cite{l}). This interplay and its applications will be discussed in detail in \cite{ad}.

As was mentioned above, the spectral theory is well understood for the primary example, the Fibonacci case, and more generally, for all Sturmian potentials. Thus it is natural to consider generalizations of Sturmian potentials and to explore whether the combinatorial approach continues to be applicable. There are a number of natural candidates:

\begin{itemize}

\item quasi-Sturmian sequences,

\item sequences obtained by codings of rotations,

\item Arnoux-Rauzy sequences.

\end{itemize}
Quasi-Sturmian sequences are essentially given by morphic images of Sturmian sequences and the corresponding Schr\"odinger operators were studied in \cite{dl4}, confirming all of the above
points. Sturmian sequences have a geometric realization as a coding of an irrational rotation on the unit circle with respect to a decomposition of the circle into two half-open intervals,
where the rotation number is equal to the length of one of the intervals. By dropping the latter condition, one obtains the more general class of sequences associated to codings of rotations. The corresponding operators display purely singular continuous zero-measure Cantor spectrum in many cases \cite{ad2,dp1,hks}. Finally, Sturmian potentials can be characterized by a scarceness
of so-called special factors, that is, aside from being defined over two symbols, there is, for each length, exactly one subword with multiple extensions to the right and one subword with
multiple extensions to the left. When considering more than two symbols, this definition leads to the class of Arnoux-Rauzy sequences, originally defined and studied in \cite{ar}. For the
corresponding operators, no results have been shown yet. Thus, the spectral analysis of these operators via the combinatorial approach mentioned above is the objective of the present paper.

To this end, we shall recall the formal definition and some basic combinatorial properties of Arnoux-Rauzy sequences in Section~2 and then study the relevant combinatorial issues, namely, linear recurrence, powers, and palindromes, in Sections~3--5, respectively. Applications of the combinatorial results obtained in these sections to the corresponding Schr\"odinger operators are then presented in Section~6.

\medskip

\noindent\textit{Acknowledgments.} We thank J.~Cassaigne for useful discussions. D.~D.\ would like to express his gratitude to the Department of Mathematics at the University of North Texas at Denton for its warm hospitality and financial support through the Texas Advanced Research Program.

\section{Basic Properties of Arnoux-Rauzy Sequences}

In this section we recall some known properties of Arnoux-Rauzy sequences and subshifts. In particular, we explain the two combinatorial descriptions of such subshifts from \cite{rz} since they will be used extensively in later sections.

We begin with some definitions. Let $\A = \{ 1 , 2 , \ldots , k \}$ with $k \ge 2$. Given a sequence $x \in \A^\N$ or $\A^\Z$, we denote by $F_x(n)$ the set of all subwords of $x$ of length $n \in \N$, that is, $F_x(n) = \{ x_j \ldots x_{j+n-1} : j \in \N \mbox{ (or $\Z$)} \}$. We write $F_x = \bigcup_{n \in \N} F_x(n)$. The complexity function $f : \N \rightarrow \N$ of $x$ is defined by $f(n) = |F_x(n)|$, where $| \cdot |$ denotes cardinality. A factor ($=$ subword) $u \in F_x(n)$ of $x$ is called right-special if it has at least two extensions to the right, that is, there are $a,b \in \A$, $a \not= b$ such that $ua, ub \in F_x(n+1)$. A left-special factor is defined analogously. If a factor is both right-special and left-special, it is called bispecial. The sequence $x$ is called an Arnoux-Rauzy sequence if

\begin{itemize}

\item $x$ is uniformly recurrent (i.e., each factor of $x$ occurs with bounded gaps),

\item $f(n) = (k-1)n + 1$,

\item each $F_x(n)$ contains exactly one right-special factor $r_n$ and one left-special factor $l_n$.

\end{itemize}
It can be shown that $r_n = l_n^R$ \cite{rz}, where the reversal $u^R$ of a word $u = u_1 \ldots u_m$ is defined by $u^R = u_m \ldots u_1$. In particular, $r_1 = l_1$ and this factor is bispecial. Observe that there is a unique symbol $a \in \A$ such that $aa \in F_x(2)$ (which is given by $a = r_1 = l_1$). We shall say that $x$ is of type $a$.

If $k = 2$, this recovers the definition of a Sturmian sequence. Hence, Arnoux-Rauzy (AR for short) sequences are a natural generalization of Sturmian sequences to larger alphabets. Given an AR sequence $x$, we define the associated AR subshift $X$ by
$$
X = \{ y \in \A^\N \mbox{ (or $\A^\Z$)} : F_y(n) = F_x(n) \mbox{ for every } n \in \N\}.
$$
By definition, we have $F_y = F_x \equiv F_X$ for every $y \in X$. When we want to be more specific (about the choice of $\N$ or $\Z$), we shall refer to $X$ as a one-sided (resp., two-sided) subshift. AR sequences and subshifts were originally defined and studied by Arnoux and Rauzy in \cite{ar}.

Since $x$ is assumed to be uniformly recurrent, $X$ is minimal. Moreover, it was shown in \cite{ar} that $X$ is uniquely ergodic, that is, $X$ admits a unique shift-invariant probability measure $\nu$. Equivalently, for every factor $u \in F_X$ and every $m \in \N$ (or $m \in \Z$), the limit
$$
d(u) = \lim_{n \rightarrow \infty} \frac{1}{n} \#_u (x_m \ldots x_{m+n-1})
$$
exists, uniformly in $m$. Here, $\#_u (v)$ denotes the number of occurrences of $u$ in $v$. The number $d(u)$ is called the frequency of $u$. We have, for every $m$,

\begin{equation}\label{measfrequ}
\nu \left( \{ y \in X : y_m \ldots y_{m + |w| - 1} = u \} \right) = d(u).
\end{equation}

Two important objects associated with such a subshift are the index sequence $(i_n) \in \A^\N$ and the characteristic sequence $(c_n) \in \A^\N$ which are defined as follows: Let $\{\varepsilon = w_1,w_2,w_3,\ldots\}$ be the set of bispecial factors ordered so that $0 = |w_1| < |w_2| < |w_3| < \cdots$. For $n \in \N$, let $i_n \in \A$ be the unique symbol so that $i_n w_n$ is right-special. The characteristic sequence $(c_n)$, on the other hand, is defined to be the unique accumulation point of the set $\{l_1, l_2, l_3, \ldots \}$ of left-special factors. Note that $(c_n)$ is an element of the one-sided subshift $X$ and hence has the same factors as $x$. Consequently, it carries all the necessary information and once we find a way of constructing $(c_n)$ from $(i_n)$, we see that the index sequence completely determines $X$.

One such construction is given by the hat algorithm from \cite{rz}. Define a function
$$
H : \A^\N \rightarrow \A^\N
$$
as follows. Set
$$
\A' = \{ 1, \ldots, k, \hat{1}, \ldots, \hat{k} \}
$$
and let $\Phi$ denote the morphism
$$
\Phi : \A' \rightarrow \A, \; \Phi(\hat{a})=\Phi(a)=a \mbox{ for every } a \in \A.
$$
Clearly, $\Phi$ extends to both $(\A')^*$ and $(\A')^\N$. With each sequence $S = (s_n) \in \A^\N$, we associate a sequence $(B_n)$ of words over the alphabet $\A'$ as follows: $B_1 = \hat{s}_1$ and, for $n > 1$, $B_n$ is obtained from $B_{n-1}$ according to the following rule. If $\hat{s}_n$ does not occur in $B_{n-1}$, then
$$
B_n = B_{n-1} \hat{s}_n \Phi(B_{n-1}).
$$
Otherwise, if $\hat{s}_n$ occurs in $B_{n-1}$, then we can write $B_{n-1} = x' \hat{s}_n y'$, where $x', y'$ are words over $\A'$ (possibly empty) and $\hat{s}_n$ does not occur in $y'$. In this case we set
$$
B_n = B_{n-1} \hat{s}_n \Phi(y').
$$
The sequence $(B_n)$ converges to a unique sequence $B \in (\A')^\N$. We set
$$
H(S) = \Phi(B).
$$
Then, the following holds:

\begin{theorem}[Risley-Zamboni \cite{rz}]\label{rzthm1}
Let $X$ be an AR subshift over $\A$. Let $I = (i_n)$ be its index sequence and $C = (c_n)$ its characteristic sequence. Then every $a \in \A$ occurs in $I$ an infinite number of times and
$$
C = H(I).
$$
Conversely, if $I = (i_n)$ is a sequence over $\A$ such that every $a \in \A$ occurs infinitely often in $I$, then $H(I)$ is the characteristic sequence of an AR subshift.
\end{theorem}

The key observation is that $\{ \Phi (B_n) : n \in \N \}$ is precisely the set of all bispecial factors (see \cite{rz}). Thus, we can re-interpret the construction above on this level: Let $w$ be one of the bispecial factors. Suppose that each symbol $a \in \A$ occurs in $w$ (this holds by minimality if $w$ is long enough). Then by the hat algorithm, for each symbol $a \in \A$, there is a positive integer $m < |w|$ (depending on $a$) such that the next bispecial factor is obtained from $w$ by adjoining to the end of $w$ a suffix of $w$ of length $m$. The quantity $m$ is one of the $k$ periods of $w$. Let $p_1 , p_2 , \ldots , p_k$ denote the $k$ periods of $w$. Then we have the following formula \cite{wz}:
$$
\sum_{i = 1}^k (p_i - 1) = (k-1) |w|.
$$
We can suppose that $p_1 \geq p_2 \geq \cdots \geq p_k$. In this case it follows that

\begin{equation}\label{piinequ}
p_i > \frac{|w|}{k}, \; 1 \le i \le k-1.
\end{equation}
In fact, if $p_{k-1} \le |w|/k$, then $p_{k-1} + p_k \le |w|$, implying that $p_1 + p_2 + \cdots + p_{k-2} \geq (k-2)|w|$, which is a contradiction since each $p_i<|w|$.

Another way of constructing the characteristic sequence is given by the following result. For each $a \in \A$, define the morphism $\tau_a$ by $\tau_a(a) = a$ and $\tau_a(b) = ab$ for $b \in \A \setminus \{a\}$.

\begin{theorem}[Risley-Zamboni \cite{rz}]\label{rzthm2}
Let $X$ be an AR subshift and let $(i_n)$ be its index sequence. For each $a \in \A$, the characteristic sequence $(c_n)$ of $X$ is given by
$$
\lim_{n \rightarrow \infty} \tau_{i_1} \circ \cdots \circ \tau_{i_n} (a).
$$
\end{theorem}

That is, the characteristic sequence admits an $S$-adic representation where the underlying morphisms are given by $\{ \tau_a : a \in \A \}$ and they are iterated in an order dictated by the index sequence.

\section{Linearly Recurrent Arnoux-Rauzy Sequences}

In this section we characterize the set of AR subshifts that are linearly recurrent. Recall that a subshift $X$ is called $K$-linearly recurrent (or $K$-LR) if there is a constant $K > 0$ such that every $w \in F_X$ is contained in every $v \in F_X$ of length $K|w|$. $X$ is called linearly recurrent (or LR) if it is $K$-LR for some $K$. Linear recurrence is a concept that has been quite popular since the late nineties and it is known to have a number of nice consequences; compare \cite{dl5,du,dhs,lp,l3,l}. For example, every linearly recurrent $X$ is uniquely ergodic and $N$-power free (i.e., $F_X$ does not contain an element of the form $u^N$).

\begin{theorem}\label{lrthm}
An AR subshift $X$ over $\A$ is linearly recurrent if and only if every letter $a \in \A$ occurs in $(i_n)$ with bounded gaps.
\end{theorem}

\noindent\textit{Remark.} This is Corollary~III.9 in \cite{rz}. One direction was stated without proof and the proof of the other direction was based on \cite[Proposition~5]{du} which turned out to be incorrect \cite{du2}. 

\begin{proof}
It clearly suffices to prove the assertion for the characteristic sequence $(c_n)$ since the LR property only depends on the set of factors of a sequence and hence is an invariant of a minimal subshift.

We first prove that if some letter $\tilde{a} \in \A$ occurs in $(i_n)$ with unbounded gaps, then $(c_n)$ is not linearly recurrent. A special case of this scenario is easy to handle: If for each $n \in \N$, there is $m \in \N$ such that $i_m = \cdots = i_{m+n-1}$, then $(c_n)$ is not LR since it is not $N$-power free for any $N$ (see \cite[Corollary~III.6]{rz} or the next section). Let us therefore assume, in addition, that there is some $N \in \N$ such that for every $a \in \A$, $a^N$ does not occur in the index sequence. Fix some $K > 0$. Let $L > kNK$ (recall that $k$ is the size of the alphabet). Then there exists $n \in \N$ such that $i_{n + j} \not= \tilde{a}$, $0 \le j \le L$. We shall show that $w_{n+L}$ is a word of length $> (K+1) | w_n |$ which does not contain an occurrence of $w_n \tilde{a}$. (Recall that $w_m$ denotes the $m$-th bispecial factor.) This implies that $(c_n)$ is not $K$-LR. Since $K$ was arbitrary, $(c_n)$ is not LR. That $w_{n+L}$ does not contain $w_n \tilde{a}$ follows from the hat algorithm and the fact that $\tilde{a}$ does not occur in $i_n , \ldots , i_{n+L}$. That $w_{n+L}$ is of length $ > K |w_n|$ follows from the fact that for each $j$, in passing from $w_{n+j}$ to $w_{n+j+1}$, one adds on a suffix; all but one of these suffixes are, by \eqref{piinequ}, of length $ > |w_n|/k$. Since $(c_n)$ is $N$-power free, each window of length $N$ in the index sequence must contain at least two distinct symbols. Thus
$$
|w_{n+jN}| > |w_n| + j \frac{|w_n|}{k},
$$
and hence
$$
|w_{n+kNK}| > |w_n| + Kk \frac{|w_n|}{k} = (K+1) |w_n|.
$$

Consider now the case where each symbol $a \in \A$ occurs in the index sequence with bounded gaps. That is, there is a number $g \in \N$ such that for every $a \in \A$ and every $m \in \N$, at least one of $i_m,\ldots ,i_{m+g-1}$ equals $a$. We have to show that $C = (c_n)$ is linearly recurrent. Recall from Theorem~\ref{rzthm2} that

\begin{equation}\label{charsequ}
C = \lim_{n \rightarrow \infty} ( \tau_{i_1} \circ \cdots \circ \tau_{i_n} ) (a) \mbox { for every } a \in \A.
\end{equation}
Fix some $\tilde{a} \in \A$ and define, for $m \in \N$,

\begin{equation}\label{intermediate}
C^{(m)} = \lim_{n \rightarrow \infty} ( \tau_{i_m} \circ \cdots \circ \tau_{i_n} ) (\tilde{a}).
\end{equation}
In order to show that $C$ is linearly recurrent, we shall employ \cite[Lemma~4]{du2} which provides a sufficient condition for LR: For each $m \in \N$, let $d_m$ be the largest gap between consecutive occurrences of a word of length $2$ in $C^{(m)}$. If the set $\{ d_m : m \in \N \}$ is bounded, then $C$ is linearly recurrent. Fix some $m \in \N$ and consider the sequence $C^{(m)}$. It is an AR sequence of type $i_m$ and its factors of length $2$ are given by
$$
F_{C^{(m)}}(2) = \{ a i_m : a \in \A \} \cup \{ i_m a : a \in \A \}.
$$
The gaps between occurrences of $a i_m$ are bounded by twice the maximal length of the gaps between occurrences of $a$ in $C^{(m+1)}$, which in turn occurs with gaps bounded by $2^g$ since at least one of $i_{m+1}, \ldots , i_{m+g}$ is equal to $a$ (the corresponding substitution produces a sequence where the gaps between successive $a$'s are bounded by $2$, and then the gaps increase under subsequent substitutions at most by a factor $2$). The same argument works for words of the form $i_m a$ and hence $d_m \le 2^g$.
\end{proof}

\section{Powers in Arnoux-Rauzy Sequences}

In this section we study the occurrences of powers in a given AR subshift $X$. As was noted in \cite{rz}, if there are arbitrarily long runs in the index sequence $(i_n)$, then there are arbitrarily high powers. Here, we shall prove the converse and provide an explicit expression for the index of $X$ (the highest power occurring in $X$) in terms of the run lengths in $(i_n)$.

To this end, we shall distinguish between two types of runs in $(i_n)$, namely, open runs and closed runs. An $r$-run in $(i_n)$ is a pair $(a,l) \in \A \times \N$ such that $i_m = a$ for $l \le m \le l + r - 1$. If the value of $r$ is understood, such an $r$-run will sometimes be  simply referred to as a run. An $r$-run $(a,l)$ is called open if $i_m \not= a$ for $1 \le m \le l-1$; otherwise, it is called closed.

Recall that the (integer) index of $X$, $\ind (X) \in \N \cup \{ \infty \}$, is defined by
$$
\ind (X) = \sup \{ p \in \N : \mbox{ there exists $u \in \A^*$ such that } u^p \in F_X \}.
$$
The following result provides an explicit formula for the index in terms of the runs in the index sequence and generalizes the corresponding result in the Sturmian case (cf., e.g., \cite{b,dl6,jp,v}).

\begin{theorem}\label{powerthm}
Let $X$ be an AR subshift over $\A$ and $(i_n)$ its index sequence. If $\ind (X)$ is defined as above, then
$$
\ind (X) = \max \{ N_1 , N_2 \},
$$
where

\begin{eqnarray*}
N_1 & = & 1 + \sup \{r \in \N : (i_n) \mbox{ contains an open $r$-run} \}, \\
N_2 & = & 2 + \sup \{r \in \N : (i_n) \mbox{ contains a closed $r$-run} \}.
\end{eqnarray*}
\end{theorem}

In particular, we obtain the following corollary which is a generalization of the corresponding result in the Sturmian case, which was proved by Mignosi in \cite{m}.

\begin{coro}
An AR subshift $X$ has finite index if and only if the runs in its index sequence are uniformly bounded.
\end{coro}

We begin by proving $\ind (X) \ge \max \{ N_1 , N_2 \}$. The key observation is given in the following lemma:

\begin{lemma}\label{keyobs}
Let $a, a_1, \ldots, a_n, \alpha,\beta \in \A$ with $a_i \not= a$, $1 \le i \le n$ and $\alpha \not= \beta$. Then $\tau_a \circ \tau_{a_1} \circ \cdots \circ \tau_{a_n} (a)$ is a prefix of $\tau_a \circ \tau_{a_1} \circ \cdots \circ \tau_{a_n} (\alpha \beta )$.
\end{lemma}

\begin{proof}
For $n = 1$, we have
$$
\tau_a \circ \tau_{a_1} (a) = a a_1 a
$$
and
$$
\tau_a \circ \tau_{a_1} (\alpha \beta) = \tau_a (a_1 \alpha' \tau_{a_1}( \beta) ) = a a_1 a \ldots,
$$
where
$$
\alpha' = \left\{ \begin{array}{cl} \varepsilon & \mbox{if } a_1 = \alpha, \\ \alpha & \mbox{if } a_1 \not= \alpha. \end{array} \right.
$$
Thus the statement is true for $n=1$. Let us now assume that the statement holds for $n$. We have

\begin{align*}
\tau_a \circ \tau_{a_1} \circ \cdots \circ \tau_{a_n} \circ \tau_{a_{n+1}}(a) & = \tau_a \circ \tau_{a_1} \circ \cdots \circ \tau_{a_n} (a_{n+1} a)\\
& = \tau_a \circ \tau_{a_1} \circ \cdots \circ \tau_{a_n} (a_{n+1}) \tau_a \circ \tau_{a_1} \circ \cdots \circ \tau_{a_n} (a).
\end{align*}
On the other hand, we have (with $\alpha', \beta'$ defined as above)

\begin{align*}
\tau_a \circ \tau_{a_1} \circ \cdots \circ \tau_{a_n} \circ \tau_{a_{n+1}}(\alpha \beta) & = \tau_a \circ \tau_{a_1} \circ \cdots \circ \tau_{a_n} (a_{n+1} \alpha' a_{n+1}  \beta')\\
& = \tau_a \circ \tau_{a_1} \circ \cdots \circ \tau_{a_n} (a_{n+1}) \tau_a \circ \tau_{a_1} \circ \cdots \circ \tau_{a_n} (\alpha' a_{n+1}  \beta').
\end{align*}
Since $\alpha \not=  \beta$, at least one of them is $\not= a_{n+1}$, say $\alpha$, and then $\alpha' = \alpha \not= a_{n+1}$. Now apply the induction hypothesis.
\end{proof}

\begin{prop}\label{run}
Suppose $a, a_1, \ldots, a_n , a_{n+1} \in \A$ with $a_i \not= a$, $1 \le i \le n+1$ and $x \in \A^\N$ has an occurrence of $a$. Then, for every $r \in \N$,
$$
\tau_a \circ \tau_{a_1} \circ \cdots \circ \tau_{a_n} \circ \tau_a^r \circ \tau_{a_{n+1}} (x) \mbox{ contains } \left( \tau_a \circ \tau_{a_1} \circ \cdots \circ \tau_{a_n} (a) \right)^{r+2}.
$$
\end{prop}

\begin{proof}
Since $x$ contains $a$ and $a_{n+1} \not= a$, $\tau_{a_{n+1}}(x)$ contains $a_{n+1} a a_{n+1}$. Thus $\tau_a^r \circ \tau_{a_{n+1}} (x)$ contains $a^r a_{n+1} a^{r+1} a_{n+1} a$. Therefore $\tau_a \circ \tau_{a_1} \circ \cdots \circ \tau_{a_n} \circ \tau_a^r \circ \tau_{a_{n+1}} (x)$ contains
$$
\left( \tau_a \circ \tau_{a_1} \circ \cdots \circ \tau_{a_n} (a) \right)^{r+1} \tau_a \circ \tau_{a_1} \circ \cdots \circ \tau_{a_n} (a_{n+1} a).
$$
By Lemma~\ref{keyobs}, $\tau_a \circ \tau_{a_1} \circ \cdots \circ \tau_{a_n} (a_{n+1} a)$ has $\tau_a \circ \tau_{a_1} \circ \cdots \circ \tau_{a_n} (a)$ as a prefix.
\end{proof}

\begin{prop}\label{powerI}
Let $X$ be an AR subshift over $\A$ and $(i_n)$ its index sequence. If $(i_n)$ contains a run $a^r$, then $X$ has a factor $u^{r+1}$. If the run $a^r$ is preceded by $a$ somewhere in $(i_n)$, then $X$ has a factor $u^{r+2}$. In particular,
$$
\ind (X) \ge \max \{ N_1 , N_2 \}.
$$
\end{prop}

\begin{proof}
Both claims follow immediately from Proposition~\ref{run} and its proof.
\end{proof}

We now aim at proving $\ind (X) \le \max \{ N_1 , N_2 \}$. This will be done by starting from a factor $u^p$, $p \ge 3$ and then performing an iterated desubstitution process which will produce an $r$-run in the index sequence. In general, we have $r = p-2$, but under certain circumstances, we have $r = p - 1$. To illustrate this procedure, let us start with an example. Suppose $X$ is an AR subshift over three symbols such that $F_X$ contains

\begin{equation}\label{red1}
(21232121232122123212123212212321)^p.
\end{equation}
Clearly, $C^{(1)} = C = (c_n)$ is of type $2$ and hence $i_1 = 2$. Thus, $C^{(1)} = \tau_2 (C^{(2)})$ and $C^{(2)}$ must contain the factor

\begin{equation}\label{red2}
(13113121311312131)^p.
\end{equation}
Now, $C^{(2)}$ is of type $1$ and hence $i_2 = 1$. Thus, $C^{(2)} = \tau_1 (C^{(3)})$ and $C^{(3)}$ must contain the factor

\begin{equation}\label{red3}
(3132313231)^{p-1} 313231323?.
\end{equation}
Observe that the last symbol in the last block cannot be desubstituted uniquely. We indicate this ambiguity by ``?'' and note that the last block is one symbol shorter than the other blocks. Next, $i_3 = 3$, $C^{(3)} = \tau_1 (C^{(4)})$ and $C^{(4)}$ must contain the factor

\begin{equation}\label{red4}
(12121)^{p-1} 1212?.
\end{equation}
The ambiguity on this level comes from the ambiguity on the previous level, that is, the ``?'' in \eqref{red3}. However, we clearly have $i_4 = 1$ and hence the last $2$ in \eqref{red4} must be followed by a $1$, so in fact $C^{(4)}$ must contain the factor

\begin{equation}\label{red4b}
(12121)^{p-1} 12121.
\end{equation}
This allows us to go up one level and replace the ``?'' in \eqref{red3} by a $1$. This deciphers all the ambiguities up to this point. Next, $C^{(4)} = \tau_1 (C^{(5)})$ and $C^{(5)}$ must contain the factor

\begin{equation}\label{red5}
(221)^{p-1} 22?
\end{equation}
and $i_5 = 2$, $C^{(5)} = \tau_2 (C^{(6)})$ and $C^{(6)}$ must contain the factor

\begin{equation}\label{red6}
(21)^{p-1} 2?.
\end{equation}
Now, $i_6$ is either $1$ or $2$, but in either case, further desubstitution yields a run of length $p-2$ in the index sequence. Namely, if $i_6 = 1$, then $i_7 = \cdots = i_{7 + (p-2) - 1} = 2$, and $i_6 = 2$ gives $i_7 = \cdots = i_{7 + (p-2) - 1} = 1$. Note that, contrary to the situation above, the ambiguities in \eqref{red5} and \eqref{red6} cannot be removed. We observe:

\begin{enumerate}

\item The desubstitution process takes $u^p$ to $a^{p-1} ?$ for some $a \in \A$ and ``?'' is either known or not. This yields at least a $(p-2)$-run in the index sequence.

\item If at no step there is an ambiguity in the desubstitution process, then $w^p$ reduces to $a^p$ and hence produces a $(p-1)$-run in the index sequence.

\end{enumerate}

These observations lead to the following lemma:

\begin{lemma}\label{twoscen}
Let $X$ be an AR subshift over $\A$ and $(i_n)$ its index sequence. Suppose there is $p \ge 3$ and a primitive $u \in \A^*$ such that $u^p \in F_X$. Then we have one of the following scenarios:

\begin{itemize}

\item[{\rm (i)}] There are $a \in \A$ and $m \in \N$ such that $i_{m+j} = a$, $1 \le j \le p-1$.

\item[{\rm (ii)}] There are $a \in \A$ and $m \in \N$ such that $i_{m+j} = a$, $1 \le j \le p-2$ and $i_j = a$ for some $j$ with $1 \le j \le m$.

\end{itemize}

\end{lemma}

\begin{proof}
Start with the word $u^p$ and perform a continued desubstitution process, as above, using $\tau_{i_1}^{-1}, \tau_{i_2}^{-1} \circ \tau_{i_1}^{-1}, \tau_{i_3}^{-1} \circ \tau_{i_2}^{-1} \circ \tau_{i_1}^{-1}, \ldots$, where the last symbol of the desubstituted word may be unknown and hence denoted by ``?''. Clearly, this process leads, after, say, $m$ steps, to a desubstituted word which has either the form $a^p$ or $a^{p-1} ?$, where $a$ is some symbol from $\A$. In particular, we must have $i_{m+j} = a$ for $1 \le j \le p-1$ (in the first case) or $1 \le j \le p-2$ (in the second case). It only remains to be shown that in the second case, we must have applied $\tau_a$ somewhere along the way. Notice that each desubstituted word results from the previous word by a deletion of a number of symbols. In particular, the word $u$ we started with must contain $a$. Consider first the case where $u$ contains at least two occurrences of $a$. Then, in order to reduce $u^2$ (the first two of the $p \ge 3$ blocks) to $a^2$, we necessarily have to apply $\tau_a$ along the way. Let us now consider the case where $u$ contains exactly one $a$. Then we do not apply $\tau_a$ until we are left with a desubstituted word of length $\le 2$. That is, either the word is $a$, in which case we are done (since the $a$ in the last block never gets deleted and hence $u^p$ reduces to $a^p$), or the word is $ab$ (or $ba$) for some $b$. In the next step, either $\tau_a$ or $\tau_b$ is applied. Remember that the last word still contains $a$ (so that it is one of $ab, ba, a?$) so that desubstitution by $\tau_b$ leads to $a^p$. That is, only desubstitution by $\tau_a$ leads to $a^{p-1} ?$.
\end{proof}

\begin{prop}\label{powerII}
Let $X$ be an AR subshift over $\A$ and $(i_n)$ its index sequence. Then
$$
\ind (X) \le \max \{ N_1 , N_2 \}.
$$
\end{prop}

\begin{proof}
This is an immediate consequence of Lemma~\ref{twoscen}. Namely, a given power $u^p \in F_X$, for some $p \ge 3$, corresponds to either an open or closed $(p-1)$-run or a closed $(p-2)$-run in the index sequence. Note that every AR subshift contains squares (e.g., $i_1 i_1$) so that powers $p < 3$ are irrelevant for the computation of the index.
\end{proof}

\begin{proof}[Proof of Theorem~\ref{powerthm}.]
The assertion follows from Propositions~\ref{powerI} and \ref{powerII}.
\end{proof}

One might also be interested in powers that occur for arbitrarily long factors. That is, define ${\rm i-ind}(X) \in \N \cup \{ \infty \}$ by
$$
{\rm i-ind}(X) = \sup \{ p \in \N : \mbox{ there exist $u_n$ with $|u_n| \rightarrow \infty$ such that } u_n^p \in F_X \}.
$$
Then, the above analysis has the following immediate consequence:

\begin{coro}\label{longpowers}
Let $X$ be an AR subshift over $\A$ and $(i_n)$ its index sequence. If ${\rm i-ind} (X)$ is defined as above, then
$$
{\rm i-ind} (X) = 2 + \limsup_{n \rightarrow \infty} e_n,
$$
where, for $n \in \N$,
$$
e_n = \max\{ l \in \N : i_{n+l-1} = i_n \}.
$$
\end{coro}

\begin{proof}
Since every symbol from $\A$ occurs in $(i_n)$, the index sequence has exactly $k$ open runs. Moreover, there is $N \in\N$ such that beyond $i_N$, there are no more open runs. In particular, for the computation of ${\rm i-ind}(X)$, only closed runs are relevant. Thus, the assertion follows in a straightforward way from Proposition~\ref{powerI}, Lemma~\ref{twoscen}, and their proofs.
\end{proof}

\section{Palindromes in Arnoux-Rauzy Sequences}

In this section we study the number of palindromes of a given length that occur in a given AR subshift $X$. Recall that a word is called a palindrome if it is the same when read backwards. Given a minimal subshift $X$, define its palindrome complexity function $p : \N \rightarrow \N_0$ by
$$
p(n) = | \{ p \in F_X : p = p^R, \; |p| = n \} |.
$$
It was shown by Droubay and Pirillo that all Sturmian subshifts have the same palindrome complexity function, namely,

\begin{equation}\label{sturmpalind}
p(n) = \left\{ \begin{array}{cl} 2 & \mbox{if $n$ is odd,} \\ 1 & \mbox{if $n$ is even,} \end{array} \right.
\end{equation}
and that Sturmian subshifts are in fact characterized by this property \cite{dp}. Here, we shall generalize the first part, namely, we show that all AR subshifts over $\A$ have the same palindrome complexity function; and we also prove that the second part does not generalize, that is, for some $k \ge 3$, there are non-AR subshifts over $\A$ that have the same palindrome complexity function as AR subshifts over $\A$.

\begin{theorem}\label{palinthm}
The palindrome complexity function $p : \N \rightarrow \N_0$ of an AR subshift $X$ over $\A$ is given by

\begin{equation}\label{arpalind}
p(n) = \left\{ \begin{array}{cl} k & \mbox{if $n$ is odd,} \\ 1 & \mbox{if $n$ is even.} \end{array} \right.
\end{equation}
\end{theorem}

\begin{proof}
We shall prove the statement

\begin{equation}\label{palstate}
\forall n \in \N : p(2n-1) = k, \; p(2n) = 1,
\end{equation}
which is equivalent to the assertion, by induction on $n$.

The case $n = 1$ is readily checked. In fact, $p(1) = k$ is obvious, and if $X$ is of type $a \in \A$, then $aa$ is the unique palindrome of length $2$ which occurs in $X$.

Now assume that \eqref{palstate} holds for $n$. Let us show \eqref{palstate} for $n+1$ by proving that if $p$ is a palindrome occurring in $X$, then $p$ admits a unique extension $apa$ to a palindrome of length $|p| + 2$.

Fix a palindrome $p \in F_X$. We show below that

\begin{equation}\label{bispecialob}
p \in F_X \mbox{ bispecial } \Rightarrow \mbox{there exists a unique $a \in \A$ such that } apa \in F_X.
\end{equation}
Now, either there exists a unique $a \in \A$ such that $apa \in F_X$, or else $p$ is bispecial (and so by \eqref{bispecialob} there exists a unique $a \in \A$ such that $apa \in F_X$); hence in either case there exists a unique $a \in \A$ such that $apa \in F_X$.

Let us show \eqref{bispecialob}. Let $a \in \A$ be the unique letter for which $ap$ is right-special (and, equivalently, $pa$ is left-special). Then $apa \in F_X$. Consider any letter $b \not= a$. Then, we have that $bp$ is not right-special, and $bpa \in F_X$ (since $pa$ is left-special), so $bpb \not\in F_X$.
\end{proof}

As we mentioned above, every minimal subshift with palindrome complexity given by \eqref{sturmpalind} is necessarily Sturmian. We are now going to show that this does not extend to the AR case, that is, there are non-AR subshifts with palindrome complexity given by \eqref{arpalind}. To this end, we consider subshifts $X_{{\rm 3iet}}$, defined over $\mathcal{A}_3$, associated with three-interval exchange transformations. These dynamical systems have the following combinatorial description, as shown by Ferenczi et al.\ \cite{fhz2}:

\begin{itemize}

\item $F_{X_{{\rm 3iet}}} (2) = \{ 12, 13, 21, 22, 31\}$.

\item If $u \in F_{X_{{\rm 3iet}}}$, then $u^R \in F_{X_{{\rm 3iet}}}$.

\item For every $n \in \N$, there are exactly two left-special words in $F_{X_{{\rm 3iet}}} (n)$, one beginning in $1$ and one beginning in $2$.

\item If $w$ is a bispecial word ending in $1$ and $w \not= w^R$, then $w2$ is left-special if and only if $w^R 1$ is left-special.

\end{itemize}

Clearly, no such subshift is an AR subshift. We have the following result:

\begin{prop}
The palindrome complexity function $p$ of $X_{{\rm 3iet}}$ is given by 

\begin{equation}\label{ar3palind}
p(n) = \left\{ \begin{array}{cl} 3 & \mbox{if $n$ is odd,} \\ 1 & \mbox{if $n$ is even.} \end{array} \right.
\end{equation}
\end{prop}

\begin{proof}
The proof is similar to the proof of Theorem~\ref{palinthm}. It follows from the proof of Proposition~2.6 in \cite{fhz2} that if $u$ is a bispecial palindrome factor, then there exists a unique symbol $a \in \mathcal{A}_3$ such that $aua$ is a factor. In fact, if $u$ begins in $1$, then $a\in \{2,3\}$, while if $u$ begins in $2$, then $a\in \{1,2\}$. So now suppose $u$ is a palindrome factor of length $n.$ We claim there exists a unique symbol $a \in \mathcal{A}_3$ such that $aua$ is a factor. If no such $a$ exists, then there exist distinct symbols $b,c \in \mathcal{A}_3$ such that $buc$ is a factor. This implies that $u$ is bispecial, so from the above there must exist $a \in \mathcal{A}_3$ such that $aua$ is a factor. Next, suppose there exist distinct symbols $b,c \in \mathcal{A}_3$ such that $bub$ and $cuc$ are factors. Then again $u$ is bispecial, and hence this cannot happen. Thus $p(n) = p(n+2)$. Since $p(1)=3$ and $p(2)=1$, the assertion follows.
\end{proof}

On the other hand, $X_{{\rm 3iet}}$ has the same factor complexity function as an AR subshift over three symbols, so one may ask whether \eqref{ar3palind} implies that

\begin{equation}\label{ar3comp}
f(n) = 2n + 1.
\end{equation}
This is not true, at least on the level of individual sequences, as demonstrated by the following result (the example is due to J.~Cassaigne \cite{c}).

\begin{prop}
There exists a sequence over $\mathcal{A}_3$ whose palindrome complexity function is given by \eqref{ar3palind}, but whose factor complexity function is not given by \eqref{ar3comp}.
\end{prop}

\begin{proof}
Let
$$
w = 121312141213121\ldots
$$
be the fixed point of the infinite substitution
$$1\mapsto 12, \; 2\mapsto 13, \; 3\mapsto 14, \ldots .
$$
Let $w'$ be the morphic image of $w$ under the map $\Theta$ where $\Theta(i) = 1^i 2 1^i 3 1^i$.
So
$$w' = 1213111211311121311112111311112131\ldots .
$$
Note that

\begin{equation}\label{observe}
22, 33, 2 1^n 2, 3 1^n 3, 3 1^n 2 1^n 3, 2 1^n 3 1^n 2 \mbox{ are not factors of } w'.
\end{equation}
We therefore have for $w'$,

\begin{equation}\label{ex1}
p(n)=1 \mbox{ for $n$ even},
\end{equation}
since by \eqref{observe} the only even length palindromes are $1^n$,

\begin{equation}\label{ex2}
p(n)=3 \mbox{ for $n$ odd},
\end{equation}
since by \eqref{observe} the only odd length palindromes are $1^n, 1^{\frac {n-1}{2}}21^{\frac {n-1}{2}}, 1^{\frac {n-1}{2}}31^{\frac {n-1}{2}}$, and

\begin{equation}\label{ex3}
f \mbox{ is not given by \eqref{ar3comp}.}
\end{equation}
For example, $f(3) = 9$. By \eqref{ex1}--\eqref{ex3}, we have palindrome complexity as in \eqref{ar3palind}, but factor complexity different from \eqref{ar3comp}.
\end{proof}

Note, however, that the sequence $w'$ above is not uniformly recurrent and hence does not induce a minimal subshift. We consider it an interesting open problem to determine all minimal subshifts $X$ over $\A$, with $k \ge 3$ arbitrary, whose palindrome complexity function is given by \eqref{arpalind}.

\section{Applications to Schr\"odinger Operators}

In this section, we discuss applications of our combinatorial results, Theorems~\ref{lrthm} and \ref{palinthm} and Corollary~\ref{longpowers}, to the spectral theory of Schr\"odinger operators.

A discrete one-dimensional Schr\"odinger operator acts in the Hilbert space $\mathcal{H} = \ell^2(\Z)$. If $\phi \in \mathcal{H}$, then $H \phi$ is given by
$$
(H\phi)(n) = \phi(n+1) + \phi(n-1) + V(n) \phi(n),
$$
where $V : \Z \rightarrow \R$. The map $V$ is called the potential. For our purposes, we can assume $V$ bounded. Then $H$ is a bounded, self-adjoint operator. Denote the spectrum of $H$ by $\sigma(H)$. Given an initial state $\phi \in \mathcal{H}$, the Schr\"odinger time evolution is given by $\phi(t) = \exp (-itH) \phi$, where $\exp (-itH)$ is given by the spectral theorem. One is interested in the question whether $\phi(t)$ will spread out in space, and if so, how fast. One possible way to tackle this issue is to study the spectral measure $\mu_\phi$ associated with $\phi$, which is defined by
$$
\langle \phi , (H - z)^{-1} \phi \rangle = \int_\R \frac{d\mu_\phi(x)}{x - z} \mbox{ for every $z$ with } {\rm Im} \, z > 0.
$$
Roughly speaking, the more continuous $\mu_\phi$, the faster the spreading of $\phi(t)$; compare, for example, \cite{bgt1,g3,l2}. Denote

\begin{align*}
\mathcal{H}_{{\rm ac}} & = \{ \phi \in \mathcal{H} : \mu_\phi \mbox{ is absolutely continuous} \}\\
\mathcal{H}_{{\rm sc}} & = \{ \phi \in \mathcal{H} : \mu_\phi \mbox{ is singular continuous} \}\\
\mathcal{H}_{{\rm pp}} & = \{ \phi \in \mathcal{H} : \mu_\phi \mbox{ is pure point} \}
\end{align*}
and
$$
\sigma_\varepsilon (H) = \sigma \left( H |_{\mathcal{H}_{\varepsilon}} \right) \mbox{ for } \varepsilon \in \{ \mbox{ac,sc,pp} \}.
$$
We'll say that $H$ has purely absolutely continuous spectrum if both $\sigma_{{\rm sc}}(H)$ and $\sigma_{{\rm pp}}(H)$ are empty, etc.

As was mentioned in the introduction, there has been a considerable amount of research dealing with the spectral properties of $H$ if $V$ displays long-range order. The totally ordered case (i.e., $V$ periodic) is well-understood \cite{teschl}. In this case, $H$ has purely absolutely continuous spectrum. If $V$ takes on only finitely many values, one popular measure for long-range order is given by the complexity function. It has then been the goal to determine the spectral properties for aperiodic potentials of low combinatorial complexity. A complete understanding has been obtained for Sturmian potentials \cite{bist,dkl} and quasi-Sturmian potentials \cite{dl4}. It turned out that in all these cases, one has purely singular continuous spectrum, supported on a Cantor set of Lebesgue measure zero. Here, a Cantor set is a closed, perfect, nowhere dense set. It is natural to conjecture that these properties are shared by other low-complexity potentials. In fact, these questions can be studied from a purely combinatorial perspective. That is, there are results that deduce singular continuous, zero-measure spectrum from purely combinatorial properties of the potential. Here we study the case of Arnoux-Rauzy potentials which provide a natural class of low-complexity potentials.

Fix a two-sided AR subshift $X$ over $\A$ with index sequence $(i_n)$ and a non-constant function $f : \A \rightarrow \R$. Denote the unique ergodic measure on $X$ by $\nu$. Each element $x$ of $X$ induces a potential via $V_x(n) = f(x_n)$. The Schr\"odinger operator with potential $V_x$ will be denoted by $H_x$. Since $X$ is minimal, we have that the spectrum and the absolutely continuous spectrum are invariants of $X$, that is, there are sets $\Sigma, \Sigma_{{\rm ac}} \subseteq \R$ such that $\sigma (H) = \Sigma$ and $\sigma_{{\rm ac}} (H) = \Sigma_{{\rm ac}}$ for every $x \in X$. The result for the spectrum follows from strong convergence and is folklore. The result on the absolutely continuous is much deeper and more recent \cite{ls}. In fact, aperiodicity implies that $\Sigma_{{\rm ac}}$ is empty \cite{k2}. Thus, to establish the desired picture, we have to show that $\Sigma$ has Lebesgue measure zero and $\sigma_{{\rm pp}}(H)$ is often/always empty.

We first turn to the zero-measure property. It is a result of Lenz that linear recurrence provides a sufficient condition:

\begin{theorem}[Lenz \cite{l}]
If $X$ is a linearly recurrent subshift and $X$ and $f$ are such that the resulting potentials $V_x$ are aperiodic, then $\Sigma$ has Lebesgue measure zero.
\end{theorem}

Combining this with our Theorem~\ref{lrthm}, we immediately obtain the following (the Cantor set properties follow from the zero-measure property by general principles):

\begin{coro}
If every letter $a \in \A$ occurs in $(i_n)$ with bounded gaps, then $\sigma(H_x)$ is a Cantor set of zero Lebesgue measure for every $x \in X$.
\end{coro}

Let us now discuss the absence of point spectrum. Both palindromes and powers allow one to prove this property. The palindrome criterion is easy to verify, but it has the slight disadvantage that it only gives generic absence of eigenvalues:

\begin{theorem}[Hof et al.\ \cite{hks}]\label{pp}
If $X$ is a minimal subshift and its palindrome complexity function obeys $\limsup_{n \rightarrow \infty} p(n) > 0$, then for a dense $G_\delta$-set of $x \in X$, we have $\sigma_{{\rm pp}}(H_x) = \emptyset$.
\end{theorem}

We immediately deduce from this and Theorem~\ref{palinthm}:

\begin{coro}
For a dense $G_\delta$-set of $x \in X$, we have $\sigma_{{\rm pp}}(H_x) = \emptyset$.
\end{coro}

On the other hand, the criterion for empty point spectrum which is based on powers is slightly more complicated to state, requires more effort to be verified, but yields a stronger conclusion. Define the set $X_n$ of elements of $X$, which have cubes of length $3n$, suitably centered around the origin, by
$$
X_n = \{ x \in X : x_{-n+j} = x_j = x_{n+j}, \; 1 \le j \le n \}.
$$
Then we have the following result (the proof is based on a Gordon-type argument \cite{g}; see, e.g., \cite{d2,dp1}):

\begin{theorem}\label{powerpp}
Suppose $\limsup_{n \rightarrow \infty} \nu (X_n) > 0$. Then, for $\nu$-almost every $x \in X$, we have $\sigma_{{\rm pp}}(H_x) = \emptyset$.
\end{theorem}

We can use this theorem and our Corollary~\ref{longpowers} to show:

\begin{coro}\label{36}
If the index sequence $(i_n)$ contains infinitely many $2$-runs, we have $\sigma_{{\rm pp}}(H_x) = \emptyset$ for $\nu$-almost every $x \in X$.
\end{coro}

\begin{proof}
Corollary~\ref{longpowers} shows that if the index sequence $(i_n)$ contains infinitely many $2$-runs, then $F_X$ contains arbitrarily long fourth powers. That is, there are $u_n \in \A^*$ with $|u_n| \rightarrow \infty$ and $u_n^4 \in F_X$. Since $X$ is aperiodic, either $u_n^4$ is right-special or one of its conjugates is right-special. Thus, we can assume without loss of generality that $u_n^4$ is right-special. It was shown in \cite[Lemma~2.2]{wz} that among all factors of length $4 |u_n|$, the right-special factor has the largest frequency. Since there are $4 (k-1) |u_n| + 1$ words of length $4 |u_n|$ whose frequencies add up to one, we infer that
$$
d(u_n^4) \ge \frac{1}{4(k-1)|u_n| + 1}.
$$
This yields, using \eqref{measfrequ},
$$
\limsup_{n \rightarrow \infty} \nu (X_n) \ge \limsup_{n \rightarrow \infty} \nu \left(X_{|u_n|}\right) \ge \limsup_{n \rightarrow \infty} \frac{|u_n|}{4(k-1)|u_n| + 1} = \frac{1}{4(k-1)} > 0.
$$
Thus, the assertion follows from Theorem~\ref{powerpp}.
\end{proof}

Corollary~\ref{36} does not cover the prominent case of the Tribonacci subshift $X_{{\rm Trib}}$, which is defined over three symbols and corresponds to the index sequence
$$
(i_n) = 1,2,3,1,2,3,1,2,3,\ldots .
$$
We shall nevertheless show that the conclusion of Corollary~\ref{36} holds for this case. By Theorem~\ref{rzthm2}, the characteristic sequence $C = (c_n)$ is given by
$$
C = \lim_{n \rightarrow \infty} \left( \tau_1 \circ \tau_2 \circ \tau_3 \right)^n (1).
$$
The substitution $S = \tau_1 \circ \tau_2 \circ \tau_3$ on $\mathcal{A}_3$ is given by
$$
S(1) = 1213121, \; S(2) = 121312, \; S(3) = 1213.
$$
Note that $S$ is primitive (i.e., there is $l \in \N$, namely $l = 1$, such that for every $a \in \mathcal{A}_3$, $S^l(a)$ contains all symbols from $\mathcal{A}_3$). Recall that a fractional power $w^q$ is a word $w^p w'$ with $p \in \N$, $w'$ a prefix of $w$, and $q = p + |w'|/|w|$. We have the following result for subshifts generated by primitive substitutions:

\begin{theorem}[Damanik~\cite{d3}]\label{fract}
Suppose the subshift $X$ is generated by a primitive substitution $S$ and $F_X$ contains a fractional power $w^q$ with $q > 3$. Then we have $\sigma_{{\rm pp}}(H_x) = \emptyset$ for $\nu$-almost every $x \in X$.
\end{theorem}

This allows us to prove the following:

\begin{coro}\label{38}
For the Tribonacci subshift $X_{{\rm Trib}}$, we have $\sigma_{{\rm pp}}(H_x) = \emptyset$ for $\nu$-almost every $x \in X_{{\rm Trib}}$.
\end{coro}

\begin{proof}
As we have seen above, $C$ is the unique fixed point of $S$ in $\mathcal{A}_3^\N$ and we have
$$
F_{X_{{\rm Trib}}} = \bigcup_{n \in \N} F_{S^n(1)}.
$$
Thus it suffices to find some $S^n(1)$ which contains $w^q$ with $q > 3$. The claim then follows from Theorem~\ref{fract}. First, $S^2(1)$ contains the word $1121$. Thus, $S^3(1)$ contains the word
$$
1213121 \, 1213121 \, 121312 \, 1213121 = ( 1213121 )^3 2 \ldots,
$$
and hence $S^4(1)$ contains
$$
( 1213121 \, 121312 \, 1213121 \, 1213 \, 1213121 \, 121312 \, 1213121 )^3 121312 \ldots,
$$
which yields a fractional power $w^q$ with $q = 3 + 3/22 > 3$.
\end{proof}

\end{document}